\documentclass{article}
\usepackage{amsfonts}
\usepackage{amsthm}
\usepackage{amssymb}
\usepackage{amsmath}
\usepackage{times}
\begin{document}

\author{Diana Dziewa-Dawidczyk\footnote{Faculty of Mathematics
and Information Science,  Warsaw University of Technology, Pl.
Politechniki 1, 00-661 Warsaw, POLAND, e-mail: ddziewa@mini.pw.edu.pl}  %
 \and Zbigniew Pasternak-Winiarski\footnote{Faculty of Mathematics
and Information Science,  Warsaw University of Technology, Pl.
Politechniki 1, 00-661 Warsaw, POLAND, Institute of Mathematics,
University of Bia\l ystok Akademicka 2, 15-267 Bia\l ystok,
Poland, e-mail: Z.Pasternak-Winiarski@mini.pw.edu.pl }}
\title{On differential completions and compactifications of a differential space}

\maketitle

\begin{abstract}
Differential completions and compactifications of differential
spaces are introduced and investigated. The existence of the
maximal differential completion and the maximal differential
compactification is proved. A sufficient condition for the
existence of a complete uniform differential structure on a given
differential space is given.
\bigskip


\bigskip
\end{abstract}

\section {Introduction}
This article is the third of the series of papers concerning
integration of differential forms and densities on differential
spaces (the first two are \cite{DPW2} and \cite{DPW1}). We
describe differential completions and differential
compactifications of differential spaces which are used in our
theory of integration.

Section 2 of the paper contains basic definitions and the
description of preliminary facts concerning theory of differential
spaces. In Section 3 we give basic definitions and describe the
standard facts concerning theory of uniform spaces. We introduce
the notion of a differential completion of a differential space.
We construct differential completions of a differential space
using families of generators of its differential structure
(Proposition 3.7, Definition 3.12). Section 4 is devoted to the
investigation of properties of differential completions. We define
some natural order in the set of all differential completions of a
given differential space. We prove that for any differential space
$(M,{\cal C})$ there exists the maximal differential completion
with respect to this order (Theorem 4.1). If for the uniform
structure defined by some family of generators the space $M$ is
complete then the appropriate differential completion of $(M,{\cal
C})$ is maximal and coincides with $(M,{\cal C})$ (Theorem 4.2).
At the end we prove that if a differential structure ${\cal C}$
posses a countable family of generators then it coincides with its
maximal differential completion (Theorem 4.3). As a corollary we
obtain general topological result about the existence on a given
topological space a uniform structure defining the initial
topology (Corollary 4.1). In Section 5 we introduce and
investigate the notion of a differential compactification of a
differential space. Similarly as in Section 4 we prove the
existence of the maximal differential compactification of a given
differential space with respect to the suitable order.

 Without any other explanation we use the following symbols:
${\mathbb N}$-the set of natural numbers; ${\mathbb R}$-the set of
reals.

\section {Differential spaces}

Let $M$ be a nonempty set and let $\mathcal{C}$ be a family of
real valued functions on $M$. Denote by $\tau_{\mathcal{C}}$ the
weakest topology on $M$ with respect to which all functions of
$\mathcal{C}$ are continuous.

A base of the topology $\tau_{\mathcal C}$ consists of sets:\\
\begin{displaymath}
(\alpha_1,\ldots,\alpha_n)^{-1}(P)=\bigcap_{i=1}^{n} \{m\in M:
a_i<\alpha_i(m)<b_i\},
\end{displaymath}
where $n\in{{\mathbb N}}$,
$a_1,\ldots,a_n,b_1,\ldots,b_n\in{\mathbb R}$, $a_i<b_i$,
$\alpha_1, \ldots,\alpha_n\in{\mathcal C}$,
$P=\{(x_1,\ldots,x_n)\in{\mathbb R}^{n};a_i<x_i<b_i,i=1,
\ldots,n\}$.\\

{\sc Definition 2.1} A function $f: M\rightarrow{\mathbb R}$ is
called \emph{a local ${\mathcal C}$-function on} $M$ if for every
$m\in M$ there is a neighborhood $V$ of $m$ and $\alpha\in
{\mathcal C}$ such that $f_{|V}=\alpha_{|V}$. The set of all local
${\mathcal C}$-functions on $M$ is denoted by ${\mathcal C}_M$.\\

Note that any function $f\in{\mathcal C}_M$ is continuous with
respect to the topology $\tau_{\mathcal C}$. Then $\tau_{{\mathcal C}_M}=
\tau_{\mathcal C}$ (see  \cite{DPW2}, \cite{DPW1}).\\

{\sc Definition 2.2} A function $f: M\rightarrow{\mathbb R}$ is
called \emph{${\mathcal C}$-smooth function on} $M$ if there exist
$n\in{\mathbb N}$, $\omega \in C^{\infty} ({\mathbb R}^{n})$ and
$\alpha_1,\ldots,\alpha_n\in{\mathcal C}$ such that
$$f=\omega \circ(\alpha_1,\ldots,\alpha_n).$$
The set of all ${\mathcal C}$-smooth functions on $M$ is denoted
by ${\rm sc}{\mathcal C}$.\\

Since ${\mathcal C}\subset sc{\mathcal C}$ and any superposition
$\omega \circ(\alpha_1,\ldots,\alpha_n)$ is
continuous with respect to $\tau_{\mathcal C}$ we obtain
$\tau_{sc{\mathcal C}}=\tau_{\mathcal C}$ (see \cite{DPW2}, \cite{DPW1}).\\

{\sc Definition 2.3} A set $\mathcal{C}$ of real functions on $M$
is said to be a \emph{(Sikorski's) differential structure} if: (i)
$\mathcal{C}$ is \emph{closed with respect to localization} i.e.
${\mathcal C} = {\mathcal C}_M$; (ii) ${\mathcal C}$ is {\sl
closed with respect to superposition with smooth functions} i.e.
${\mathcal C}={\rm sc}{\mathcal C}$.\\

In this case a pair $(M,{\mathcal C})$ is said to be a \emph{
(Sikorski's) differential space} (see \cite{sik}). Any element of
${\mathcal C}$ is called a \emph{smooth function on} $M$ (with respect to
${\mathcal C}$).\\

{\sc Proposition 2.1.} {\it The intersection of any family of
differential structures defined on a set $M\ne\emptyset$ is a
differential structure on $M$.}

\smallskip
For the proof see \cite{DPW2}, \cite{DPW1}, Proposition 2.1.\hfill $\Box$\\

Let ${\mathcal F}$ be a set of real functions on $M$. Then, by
Proposition 2.1, the intersection ${\mathcal C}$ of all
differential structures on $M$ containing ${\mathcal F}$ is a
differential structure on $M$. It is the smallest differential
structure on $M$ containing ${\mathcal F}$. One can easy prove
that ${\mathcal C}=(sc{\mathcal F})_M$ (see \cite{wal}). This
structure is called \emph{the differential structure generated by}
${\mathcal F}$ and is denoted by $gen({\cal F})$. Functions of
${\mathcal F}$ are called \emph{generators} of the differential
structure ${\mathcal C}$. We have also $\tau_{(sc{\mathcal
F})_M}=\tau_{sc{\mathcal F}}=\tau_{\mathcal F}$
(see remarks after Definitions 2.1 and 2.2).\\

\smallskip
Let $(M,{\mathcal C})$ and $(N, {\mathcal D})$ be differential
spaces. A map $F: M \rightarrow N$ is said to be \emph{smooth} if
for any $\beta\in{\mathcal D}$ the superposition $\beta\circ
F\in{\mathcal C}$. We will denote the fact that $\mathcal F$ is
smooth writing
$$F:(M,{\mathcal C})\to(N,{\mathcal D}).$$
If $F:(M,{\mathcal C})\to(N,{\mathcal D})$ is a bijection and
$F^{-1}:(N,{\mathcal D})\to(M,{\mathcal C})$ then $F$ is called
\emph{a diffeomorphism}.

\smallskip
If $A$ is a nonempty subset of $M$ and ${\mathcal C}$ is a differential
structure on $M$ then ${\mathcal C}_A$ denotes the differential structure
on $A$ generated by the family of restrictions
$\{\alpha_{|A}:\alpha\in{\mathcal C}\}$. The differential space
$(A,{\mathcal C}_A)$ is called \emph{a differential subspace} of
$(M,{\mathcal C})$. One can easy prove the following\\

{\sc Proposition 2.2.} {\it Let $(M,{\mathcal C})$ and $(N, {\mathcal D})$
be differential spaces and let $F:M\to N$. Then $F:(M,{\mathcal C})\to
(N,{\mathcal D})$ iff} $F:(M,{\mathcal C})\to(F(M),F(M)_{\mathcal D}).$

\medskip
If the map $F:(M,{\mathcal C})\to(F(M),F(M)_{\mathcal D})$ is a
diffeomorphism then we say that $F:M\to N$ is \emph{a
diffeomorphism onto its range} (in $(N,{\mathcal D}))$. In
particular the natural embedding $A\ni m\mapsto i(m):= m\in M$ is
a diffeomorphism of $(A,{\mathcal C}_A)$ onto its range in
$(M,{\mathcal C})$.

\bigskip
If $\{(M_i,{\mathcal C}_i)\}_{i\in I}$ is an arbitrary family of
differential spaces then we consider the Cartesian product
$\prod\limits_{i\in I}M_i$ as a differential space with the
differential structure $\hat{\bigotimes\limits_{i\in I}}{\mathcal
C}_i$ generated by the family of functions ${\mathcal
F}:=\{\alpha_i\circ pr_i:i\in I,\alpha_i\in{\mathcal C}_i\}$,
where $\prod\limits_{i\in I}M_i\ni(m_i)\mapsto pr_j((m_i))=:m_j\in
M_j$ for any $j\in I$. The topology
$\tau_{\hat{\bigotimes\limits_{i\in I}}{\mathcal C}_i}$ coincides
with the standard product topology on $\prod\limits_{i\in I}M_i$.
We will denote the differential structure
$\hat{\bigotimes\limits_{i\in I}}C^\infty({\mathbb R})$ on
${\mathbb R}^I$ by $C^\infty({\mathbb R}^I)$. In the case when $I$
is an $n$-element finite set the differential structure
$C^\infty({\mathbb R}^I)$ coincides with the ordinary differential
structure $C^\infty({\mathbb R}^n)$ of all real-valued functions
on ${\mathbb R}^n$ which posses partial derivatives of any order
(see \cite{sik}). In any case a function $\alpha:{\mathbb
R}^I\to{\mathbb R}$ is an element of $C^\infty({\mathbb R}^I)$ iff
for any $a=(a_i)\in{\mathbb R}^I$ there are $n\in{\mathbb N}$,
elements $i_1,i_2,\ldots,i_n\in I$, a set $U$ open in ${\mathbb
R}^n$ and a function $\omega\in C^\infty({\mathbb R}^n)$ such that
$a\in U[i_1,i_2,\ldots,i_n]:=\{(x_i)\in{\mathbb R}^I:
(x_{i_1},x_{i_2},\ldots,x_{i_n})\in U\}$ and for any $x=(x_i)\in
U[i_1,i_2,\ldots,i_n]$ we have
$$\alpha(x)=\omega(x_{i_1},x_{i_2},\ldots,x_{i_n}).$$

\bigskip
Let ${\mathcal F}$ be a family of generators of a differential
structure ${\mathcal C}$ on a set $M$. \emph{The generator
embedding} of the differential space $(M,{\mathcal C})$ into the
Cartesian space defined by ${\mathcal F}$ is a mapping
$\phi_{{\mathcal F}}: (M,{\mathcal C}) \rightarrow ({\mathbb
R}^{{\mathcal F}},C^{\infty} ({\mathbb R}^{{\mathcal F}}))$ given
by the formula
$$\phi_{{\mathcal F}}(m)=(\alpha(m))_{\alpha \in {\mathcal F}}$$
(for example if ${\mathcal F} = \{\alpha_1,\alpha_2,\alpha_3\}$
then $\phi_{{\mathcal
F}}(m)=(\alpha_1(m),\alpha_2(m),\alpha_3(m))\in {\mathbb R}^3
\cong{\mathbb R}^{\mathcal F}$). If ${\mathcal F}$ separates
points of $M$ the generator embedding is a diffeomorphism onto its
image. On that image we consider a differential structure of a
subspace of $({\mathbb R}^{\mathcal F},C^{\infty}
({\mathbb R}^{{\mathcal F}}))$ (see \cite{DPW1}, Proposition 2.3).\\

\bigskip
Let $M$ be a group (a ring, a field, a vector space over the field
${\mathbb K}$). A differential structure ${\mathcal C}$ on $M$ is
said to be \emph{a group} (\emph{ring, field, vector space})
\emph{differential structure} if the suitable group (ring, field,
vector space) operations are smooth with respect to ${\mathcal
C}$, ${\mathcal C}\otimes{\mathcal C}$ and ${\mathcal C}_{\mathbb
K}$, where ${\mathcal C}_{\mathbb K}$ is a field differential
structure on ${\mathbb K}$. In this case the pair $(M,{\mathcal
C})$ is called \emph{a differential group} (\emph{ring field,
vector space}). If ${\mathbb K}={\mathbb R}$ or ${\mathbb
K}={\mathbb C}$ we take ${\mathcal C}_{\mathbb
K}=C^\infty({\mathbb K})$ as a standard field differential
structure. \\

{\sc Proposition 2.3.} {\it Let $V$ be a vector space over
${\mathbb R}$  and let ${\mathcal F}$ be a family of constant
functions and linear functionals defined on $V$. Then the
differential structure ${\mathcal C}$ generated by ${\mathcal F}$
on $V$ is a vector space differential structure}.

\smallskip
For the proof see \cite{DPW2}, Proposition 2.3.\hfill $\Box$\\

\bigskip
{\sc Definition 2.4.} By \emph{a tangent vector} to a differential
space $(M, \mathcal{C})$ at a point $m \in M$ we call an ${\mathbb
R}$-linear mapping $v:\mathcal{C} \rightarrow {\mathbb R}$
satisfying the Leibniz condition: $v(\alpha \cdot \beta)=
\alpha(m)v(\beta) + \beta(m)v(\alpha)$ for any  $\alpha, \beta \in
\mathcal{C}$. We denote by $T_{m}M$ the set of all vectors tangent
to $(M,\mathcal{C})$ at the point $m\in M$ and call it \emph{the
tangent space to $(M,\mathcal{C})$ at the point} $m$. The union
$TM:=\bigcup\limits_{m\in M}T_mM$ is called \emph{the tangent
space to} $(M,\mathcal{C})$.

\smallskip
The set $TM$ can be endowed with a differential structure in the
following standard way. We define \emph{the projection} $\pi:TM\to
M$ such that for any $m\in M$ and any $v\in T_mM$
$$\pi(v)=m.$$
For any $\alpha\in\mathcal{C}$ we define \emph{the differential}
(or \emph{the exterior derivative}) of $\alpha$ as a map
$d\alpha:TM\to{\mathbb R}$ given by the following formula
$$d\alpha(v):=v(\alpha), \qquad v\in TM.$$
Then we define ${\mathcal TC}$ as the differential structure on
$TM$ generated by the family of functions ${\mathcal
TC}_0:=\{\alpha\circ\pi:\alpha\in{\mathcal C}\} \cup\{
d\alpha:\alpha\in{\mathcal C}\}$. From now on we will consider
$TM$ as a differential space with the differential structure
${\mathcal TC}$.

For any $m\in M$ we will denote by $d\alpha_m$ the restriction
$d\alpha_{|T_mM}$. It is clear that $d\alpha_m$ is a linear
functional on $T_mM$.

\smallskip
We have also that $\pi:(TM,{\mathcal TC})\to(M,{\mathcal C})$.
Then $\pi$ is continuous and for any $U\in\tau_{\mathcal C}$ the
set $TU:=\bigcup\limits_{m\in U} T_mM=\pi^{-1}(U)$ is open in $TM$
($TU\in\tau_{\mathcal TC}$). It can be proved that $TU$ is
(isomorphic to) a tangent space to the differential space
$(U,{\mathcal C}_U)$.

\bigskip
{\sc Theorem 2.1.} {\it If $(M,{\mathcal C})$ is a differential
space then for any $m\in M$ the pair $(T_mM,{\mathcal TC}_{T_mM})$
is a differential vector space and $T_mM$ is a Hausdorff space
(with respect to the topology induced by ${\mathcal TC}_{T_mM})$}.

\smallskip
For the proof see \cite{DPW2}, Theorem 3.1.\hfill $\Box$\\

{\sc Proposition 2.4.} {\it Let $(M,{\mathcal C})$ and
$(N,{\mathcal D})$ be differential spaces and let $F:(M,{\mathcal
C})\to(N,{\mathcal D})$. Then for any $v\in TM$ the linear
functional $TF(v):{\mathcal D}\to{\mathbb R}$ given by the formula
\begin{equation} \label{TF}
[TF(v)](\beta):=v(\beta\circ F),\qquad \beta\in{\mathcal D},
\end{equation}
is an element of $T_{F(\pi_M(v))}$, where $\pi_M:TM\to M$ is the
natural projection.

\smallskip
Proof.} We should show that $TF(v)$ fulfils Leibniz condition at
the point $F(m)$, where $m=\pi_M(v)$. For any $\beta,\gamma\in{\mathcal D}$
we have
$$[TF(v)](\beta\gamma)=v((\beta\gamma)\circ F)=v((\beta\circ F)(\gamma\circ F))$$
$$=v(\beta\circ F)\cdot(\gamma\circ F)(m)+(\beta\circ F)(m)\cdot
v(\gamma\circ F)$$
$$=[TF(v)](\beta)\cdot\gamma(F(m))+\beta(F(m))\cdot[TF(v)](\gamma).$$
\hfill $\Box$\\

{\sc Definition 2.5} Let $(M,{\mathcal C})$ and $(N,{\mathcal D})$
be differential spaces and let $F:(M,{\mathcal C})\to(N,{\mathcal
D})$. The map $TF:TM\to TN$ given in Proposition 2.4 by the
formula (\ref{TF}) is called \emph{the map tangent to} $F$.\\

{\sc Proposition 2.5.} {\it If $F:(M,{\mathcal C})\to(N,{\mathcal
D})$ then $TF:(TM,{\mathcal TC})\to(TN,{\mathcal TD})$ and
$\pi_N\circ TF=F\circ\pi_M$, where $\pi_M:TM\to M$ and
$\pi_N:TN\to N$ are natural projections.

\smallskip
Proof.} The second part of the thesis follows immediately from
Proposition 2.4 because for any $v\in TM$ the vector $TF(v)\in
T_{F(\pi_M(v))}$. Then $\pi_N(TF(v))=F(\pi_M(v))$. For the proof
of the first part of the thesis it enough to show that for any
$\kappa\in{\mathcal TD}_0$ the superposition $\kappa\circ
TF\in{\mathcal TC}$. Let us consider the case
$\kappa=\beta\circ\pi_N$. We have
$$\kappa\circ TF=\beta\circ(\pi_N\circ TF)=(\beta\circ F)\circ\pi_M.$$
Since $(\beta\circ F)\in{\mathcal C}$ we obtain $\kappa\circ
TF\in{\mathcal TC}$.

Now suppose that $\kappa=d\beta$, where $\beta\in{\mathcal D}$.
For any $v\in TM$ we have
$$\kappa\circ TF(v)=d\beta(TF(v))=[TF(v)](\beta)=v(\beta\circ F)
=d(\beta\circ F)(v).$$
Hence $\kappa\circ TF=d(\beta\circ F)\in{\mathcal TC}$. \hfill $\Box$\\

Let us consider the differential space $({\mathbb
R}^I,C^\infty({\mathbb R}^I))$. The differential structure
$C^\infty({\mathbb R}^I)$ is generated by the family of
projections ${\mathcal F}:=\{pr_i\}_{i\in I}$, where
$$pr_j((x_i)):=x_j,\qquad (x_i)\in{\mathbb R}^I,\ \ j\in I.$$
For any $x=(x_i),\ v=(v_i)\in{\mathbb R}^I$ the functional
$\vec{v}:C^\infty({\mathbb R}^I)\to{\mathbb R}$ given by the
formula
$$\vec{v}(\alpha):=\sum\limits_{i\in I}v_i\frac{\partial\alpha}{\partial x_i}(x)$$
is well defined (in some neighbourhood of $x$ the function
$\alpha$ depends on finite number of variables $x_i$) and is a
vector tangent to ${\mathbb R}^I$ at $x$. On the other hand, if
$u\in T_x{\mathbb R}^I$ and for any $i\in I$ we denote
$v_i:=u(pr_i)$ then for any $\alpha\in C^\infty({\mathbb R}^I)$ we
have $\vec{v}(\alpha)=u(\alpha)$. Then we identify the set
$T_x{\mathbb R}^I$ with $\{x\}\times{\mathbb R}^I$. Consequently
we identify the set $T{\mathbb R}^I$ with ${\mathbb
R}^I\times{\mathbb R}^I$. In this case the differential structure
${\mathcal T}C^\infty({\mathbb R}^I)$ is generated by the family
of functions ${\mathcal TF}:=\{pr_i\circ\pi\}_{i\in
I}\cup\{dpr_i\}_{i\in I}$, where
$$\pi(x,v)=x,\qquad (x,v)\in{\mathbb R}^I\times{\mathbb R}^I.$$
Hence for any $j\in I$
$$pr_j\circ\pi((x_i),(v_i))=x_j\ \ {\rm and}\ \ dpr_j((x_i),(v_i))=v_j.$$
It means that ${\mathcal T}C^\infty({\mathbb
R}^I)=C^\infty({\mathbb R}^I\times{\mathbb R}^I)$ and consequently
for any $x\in{\mathbb R}^I$ the differential structure ${\mathcal
T}C^\infty({\mathbb R}^I)_{T_x{\mathbb R}^I}$ is generated by the
family of projections $\{pr'_i:\{x\}\times{\mathbb R}^I\to{\mathbb
R}\}_I$, where
$$pr'_j(x,(v_i))=v_j.$$
Then we can identify ${\mathcal T}C^\infty({\mathbb
R}^I)_{T_x{\mathbb R}^I}$ with $C^\infty({\mathbb R}^I)$.

\medskip
Let $\phi_{{\mathcal F}}: (M,{\mathcal C}) \rightarrow ({\mathbb
R}^{{\mathcal F}},C^{\infty} ({\mathbb R}^{{\mathcal F}}))$ be the
generator embedding of the differential Hausdorff space
$(M,{\mathcal C})$ defined by some family of generators ${\mathcal
F}$. Then we can identify differential spaces $(M,{\mathcal C})$
and $(\phi_{{\mathcal F}}(M),C^{\infty} ({\mathbb R}^{\mathcal
F})_{ \phi_{{\mathcal F}}(M)})$ ($\phi_{{\mathcal F}}$ is a
diffeomorphism). We also identify tangent spaces $T_mM$ and
$T_{\phi_{\mathcal F}(m)}\phi_{\mathcal F}(M)$ using \emph{the
tangent map} $T\phi_{\mathcal F}$ (for any $\alpha\in C^{\infty}
({\mathbb R}^{\mathcal F})_{\phi_{{\mathcal F}}(M)}$).

\bigskip
{\sc Theorem 2.2.} {\it Let $I$ be an arbitrary nonempty set and
let $X$ be a nonempty subset of the Cartesian space ${\mathbb
R}^I$. Then for any $x=(x_i)\in X$ the space $T_xX$ tangent to the
differential space $(X,C^\infty({\mathbb R}^I)_X)$ at the point
$x$ is a closed subspace of the space $T_x{\mathbb R}^I$ tangent
to the differential space $({\mathbb R}^I,C^\infty({\mathbb
R}^I))$ at $x$}.

\smallskip
For the proof see \cite{DPW2}, Theorem 3.2.\hfill $\Box$\\

\bigskip
  A map $X:M\to TM$ such that for any $m\in M$ the value $X(m)\in T_mM$ is called
\emph{a vector field} on $M$. A vector field $X$ on $M$ is
\emph{smooth} if $X:(M,{\mathcal C})\to(TM,{\mathcal TC})$.

\section {Uniform structures and completions of a differential space defined
by families of generators}

For the general theory of uniform structures and completions see
\cite{eng}, Chapter 8 or \cite{Bour}. It is also described in
\cite{pas} and \cite{DPW1}. Here we start with the definition of
the uniform structure given on a differential space by a family
${\mathcal F}$ of generators of its differential structure.

Let ${\mathcal F}$ be a family of real-valued functions on a set
$M$ and let $(M,{\mathcal C})$ be a differential space such that
$\mathcal{C}=({\rm sc}{\mathcal F})_M$ and $(M,\tau_{\mathcal C})$
is a Hausdorff space (the last is true iff the family ${\mathcal
C}$ separates points in $X$ iff the family ${\mathcal F}$
separates points in $X$). On the set $M$ the family ${\mathcal F}$
defines the uniform structure ${\mathcal U}_{\mathcal F}$ such
that the base ${\mathcal B}$ of ${\mathcal U}_{\mathcal F}$
is given as follows:\\
\begin{equation} \label{US1}
{\mathcal B}=\{V(f_{1},\ldots,f_{k},\varepsilon)\subset M\times M;
k\in{\mathbb N};f_{1},\ldots,f_{k}\in{\mathcal F},
\varepsilon>0\},
\end{equation}
where
$$V(f_{1},\ldots,f_{k},\varepsilon)=\{(x,y)\in M\times M:\forall
1\leq i\leq k\quad |f_{i}(x)-f_{i}(y)|<\varepsilon\}$$ (see
\cite{DPW1}, Proposition 3.1).\\

\smallskip
{\sc Definition 3.1} The uniform structure ${\mathcal U}$ on a set
$M$ is said to be {\sl a differential uniform structure} on the
differential space $(M,{\mathcal C})$ if there exist a family
${\mathcal F}$ of generators of ${\mathcal C}$ such that
${\mathcal U}={\mathcal U}_{\mathcal F}$, where  ${\mathcal
U}_{\mathcal F}$ is defined by the base (\ref{US1}). The uniform
space $(M,{\mathcal U}_{\mathcal F})$ is said to be {\sl the
uniform space given by the family of generators} ${\mathcal F}$.

\bigskip
If we have two different families ${\mathcal F}_1$ and ${\mathcal
F}_2$ of generators of a differential space $(M,{\mathcal C})$,
then the uniform structures ${\mathcal U}_{{\mathcal F}_1}$ and
${\mathcal U}_{{\mathcal F}_1}$ can be different too.

\smallskip
{\sc Example 3.1} Let $M={\mathbb R}$, ${\mathcal
C}=C^\infty({\mathbb R})$, ${\mathcal F}_1=\{id_{\mathbb R}\}$ and
${\mathcal F}_2=\{id_{\mathbb R},f\}$, where
$$id_{\mathbb R}(x)=x,\quad{\mathrm and}\quad f(x)=x^2,\qquad x\in{\mathbb R}.$$
Then does not exists $\varepsilon >0$ such that $V(id_{\mathbb
R},\varepsilon)\subset V(f,1)$. Hence $V(f,1)\notin {\mathcal
U}_{{\mathcal F}_1}$ and ${\mathcal U}_{{\mathcal F}_1}\neq
{\mathcal U}_{{\mathcal F}_2}$. \hfill $\Box$\\

\bigskip
If ${\mathcal F}$ is a family of generators of a differential
structure ${\mathcal C}$ on a set $M$ then we define a uniform
structure ${\mathcal{U}}_{\mathcal{TF}}$ on the space $TM$ tangent
to the differential space $(M,\mathcal{C})$ using the family of
real-valued functions
$${\mathcal TF}=\{f \circ \pi : f \in{\mathcal F}\}\cup
\{ df: f \in {\mathcal F} \},$$ where $\pi : TM \rightarrow M$ is
the natural projection and $df:TM \rightarrow{\mathbb R}$, $df(v)
= v(f)$. As we know from the previous section, the family
$\mathcal{TF}$ generates the natural differential structure
$\mathcal{TC}$ on $TM$. The base $\mathcal{D}$ of
$U_{\mathcal{TF}}$ is given by:
$${\mathcal D} = \{V(\pi \circ f_{1}, \ldots , \pi \circ f_{k},
 df_{k+1}, \ldots , df_{m},\varepsilon) \subset TM \times TM:k,m
\in{\mathbb N},$$
$$f_{1},\ldots, f_{m}\in {\mathcal F},
\varepsilon>0 \}.$$

\bigskip
Let $(X,{\mathcal U})$, and $(Y,{\mathcal V})$ be uniform spaces.\\

\smallskip
{\sc Definition 3.2}  A mapping $f:X\to Y$ is said to be {\sl
uniform} with respect to uniform structures ${\mathcal U}$ and
${\mathcal V}$ if
$$\forall V\in{\mathcal V}\ \exists U\in{\mathcal U}\ \forall x,x'\in X\
[|x-x'|<U\Rightarrow |f(x)-f(x')|<V].$$
In other words for every $V\in{\mathcal V}$ there is $U\in{\mathcal U}$ such that
$U\subset (f\times f)^{-1}(V)$. We denote it by
$$f:(X,{\mathcal U})\to(Y,{\mathcal V}).$$

It is easy to prove that:\\
(i) any uniform mapping $f:(X,{\mathcal U})\to(Y,{\mathcal V})$ is
continuous with respect to topologies $\tau_{\mathcal U}$ and $\tau_{\mathcal V}$;\\
(ii) a superposition of uniform mappings is a uniform mapping.\\

\bigskip
We can give criteria of the uniformity:\\

\smallskip
{\sc Theorem 3.1} {\it Let $f:X\to Y$ and let ${\mathcal U}$ and
${\mathcal V}$ be uniform structures on $X$ and $Y$ respectively.
Then the following conditions are equivalent:\\

(a) $f:(X,{\mathcal U})\to(Y,{\mathcal V})$.\\

(b) If ${\mathcal B}$ and ${\mathcal D}$ are bases of ${\mathcal U}$ and ${\mathcal V}$
respectively then for each $V\in{\mathcal D}$ there exists $U\in{\mathcal B}$ such that
$U\subset (f\times f)^{-1}(V)$.\\

(c) For every pseudometric $\varrho$ in $Y$ uniform with respect
to ${\mathcal V}$, a pseudometric $\sigma$ in $X$ given by the
formula
$$\sigma (x,y)=\varrho(f(x),f(y)),\qquad x,y\in X,$$
is uniform with respect to the uniform structure} ${\mathcal U}$.

\smallskip
For the proof see \cite{eng}.

\bigskip
A mapping $f$, that is uniform with respect to uniform structures
${\mathcal U}$ and ${\mathcal V}$ could not be uniform with
respect to another uniform structures $\overline{{\mathcal U}}$
and $\overline{{\mathcal V}}$ defined on $X$ and $Y$ respectively
even if topologies $\tau_{\mathcal U}=\tau_{\overline{\mathcal
U}}$ and $\tau_{\mathcal
V}=\tau_{\overline{\mathcal V}}$.\\

{\sc Example 3.2} Let $M={\mathbb R}$,
$\mathcal{C}=C^{\infty}({\mathbb R})$, ${\mathcal F}_{1} =
\{id_{{\mathbb R}}\},\ {\mathcal F}_{2} = \{id_{{\mathbb R}},f\}$,
where $f(x)=x^{2}$, $x \in{\mathbb R}$. Let ${\mathcal V}$ be a
standard uniform structure on ${\mathbb R}$. Then the map $f$ is
uniform with respect to ${\mathcal U}_{{\mathcal F}_2}$ and
${\mathcal V}$, but it is not uniform with respect to ${\mathcal
U}_{{\mathcal F}_1}$ and ${\mathcal V}$. In fact, if $V
=\{(x,y)\in{\mathbb R}: |x - y|< \varepsilon \}\in{\mathcal V}$,
then does not exists $U \in {\mathcal U}_{{\mathcal F}_1}$ such
that ${\mathcal U} \subset (f\times f)^{-1}(V)$.

\bigskip
{\sc Definition 3.3} A bijective mapping $f:(X,{\mathcal
U})\to(Y,{\mathcal V})$ is {\sl a uniform homeomorphism} if
$f^{-1}$ is a uniform mapping. Then we say that $(X,{\mathcal U})$
and $(Y,{\mathcal V})$ are {\sl uniformly homeomorphic}.\\

\smallskip
By (i) it is obvious that if $f:(X,{\mathcal U})\to(Y,{\mathcal
V})$ is a uniform homeomorphism then $f$ is a homeomorphism of the
topological spaces $(X,\tau_{\mathcal U})$ and $(Y,\tau_{\mathcal
V})$.\\

\bigskip
Let $(X,{\mathcal U})$ be a uniform space and let $A\subset X$.
Then the family ${\mathcal U}_A:=\{(A\times A)\cap U:U\in{\mathcal
U}\}$ is a uniform structure on $A$. The uniform space
$(A,{\mathcal U}_A)$ is called {\sl the uniform subspace} of the
uniform space $(X,{\mathcal U})$. Note that if ${\mathcal F}$ is a
family of generators of a differential structure ${\mathcal C}$ on
a set $M$, $A\subset M$ and ${\mathcal
F}_{|A}=\{f_{|A}:f\in{\mathcal F}\}$, then the uniform space
$(A,{\mathcal U}_{{\mathcal F}_{|A}})$ is a uniform subspace of
the uniform space $(M,{\mathcal U}_{\mathcal F})$.

\bigskip
{\sc Definition 3.4} Let $(X,{\mathcal D})$ be a uniform space and
$V\in{\mathcal D}$. A set $U\subset X$ is said to be {\sl small of
rank} $V$ if $\exists x\in U \ \forall y\in U\ [(x,y)\in V]$
(see \cite{DPW1}, Definition 2.2).\\

\smallskip
If we define {\sl the ball} $K(x,V)$ as a set:
$$K(x,V)=\{y\in X:(x,y)\in V\}$$
then a set $U\subset X$ is small of rank $V$ iff $\exists x\in U\
[U\in K(x,V)]$.

If $F\subset X$ and $V\in{\mathcal D}$ we define {\sl the
V-neighbourhood} of $F$ as a set
$$K(F,V):=\bigcup\limits_{x\in F}K(x,V)=\{y\in X:\exists x\in F\ [(x,y)\in V]\}.$$

\bigskip
{\sc Definition 3.5} A nonempty family ${\mathcal F}$ of subsets
of a set $X$ is said to be {\sl a filter on} $X$ if:\\

(F1) $(F\in{\mathcal F}\ \wedge\ F\subset U\subset X)\ \Rightarrow\ (U\in{\mathcal F})$;\\

(F2) $(F_1,F_2 \in{\mathcal F})\Rightarrow(F_1 \cap F_2 \in{\mathcal F})$;\\

(F3) $\emptyset\notin{\mathcal F}$.\\

{\sc Definition 3.6} A filtering base on $X$ is a nonempty family
${\mathcal B}$ of subsets of $X$ such that\\

(FB1) $\forall A_1,A_2 \in{\mathcal B}\ \exists A_3\in{\mathcal B}\
[A_3\subset A_1\cap A_2]$;\\

(FB2) $\emptyset\notin{\mathcal B}$.\\

\smallskip
If ${\mathcal B}$ is a filtering base on $X$ then
$${\mathcal F}=\{F\subset X:\exists A\in{\mathcal B}\ [A\subset F]\}$$ is a filter
on $X$. It is called {\sl the filter defined by} ${\mathcal B}$ and ${\mathcal B}$
is called {\sl the base of} ${\mathcal F}$.\\

\smallskip
{\sc Proposition 3.1} {\it If $\{{\mathcal F}_i\}_{i\in I}$ is the
family of filters on the set $X$ then the intersection
$\bigcap\limits_{i\in I}{\mathcal F}_i$ is a filter on $X$.

\smallskip
Proof.} It is obvious that $\bigcap\limits_{i\in I}{\mathcal F}_i$
fulfils (F3). Suppose now that $F\in\bigcap\limits_{i\in
I}{\mathcal F}_i$. Then $F\in{\mathcal F}_i$ for any $i\in I$ and
if $F\subset U\subset X$ we obtain $U\in{\mathcal F}_i$. Hence
$U\in\bigcap\limits_{i\in I}{\mathcal F}_i$. It means that
$\bigcap\limits_{i\in I}{\mathcal F}_i$ fulfils (F1).

Let us consider now two arbitrary elements
$F,G\in\bigcap\limits_{i\in I}{\mathcal F}_i$. For any $i\in I$ we
have $F,G\in{\mathcal F}_i$ and therefore (by (F2)) $F\cap
G\in{\mathcal F}_i$. Hence $F\cap G\in\bigcap\limits_{i\in
I}{\mathcal F}_i$. It means that $\bigcap\limits_{i\in I}{\mathcal
F}_i$ fulfils (F2). \hfill $\Box$

\bigskip {\sc Definition 3.7} Let $X$ be a topological space. We
say that a filter $\mathcal{F}$ on $X$ {\sl is convergent to}
$x\in X$ (${\mathcal F} \rightarrow x$) if  for any neighbourhood
$U$ of $x$ there exists $F\in{\mathcal F}$ such that $F \subset
U$ (i.e. $U\in{\mathcal F}$).\\

\smallskip
{\sc Proposition 3.2} {\it If $X$ is a topological space and for
any $i\in I$ the filter ${\mathcal F}_i\rightarrow x$ then
$\bigcap\limits_{i\in I}{\mathcal F}_i\rightarrow x$.

\smallskip
Proof.} For any neighbourhood $U$ of $x$ and any $i\in I$ we have
$U\in{\mathcal F}_i$. Hence $U\in\bigcap\limits_{i\in I}{\mathcal
F}_i$. It means that $\bigcap\limits_{i\in I}{\mathcal
F}_i\rightarrow x$. \hfill $\Box$

\smallskip
{\sc Definition 3.8} Let $(X,{\mathcal U})$ be a uniform space. A
filter ${\mathcal F}$ on $X$ is {\sl a Cauchy filter} if
$$\forall V\in{\mathcal U}\ \ \exists F\in{\mathcal F}\ \ [F\times F\subset V].$$
We say that two Cauchy filters ${\mathcal F}_1$ and ${\mathcal
F}_2$ \emph{are in the relation} $R$ if
$$\forall V\in{\mathcal U}\ \ \exists F_1\in{\mathcal F}_1,
F_2\in{\mathcal F}_2\ \ [F_1\times F_2\subset V].$$

\smallskip
{\sc Proposition 3.3} {\it Two filters ${\mathcal F}_1$ and
${\mathcal F}_2$ on the uniform space $(X,{\mathcal U})$ are in
the relation $R$ iff ${\mathcal F}_1$, ${\mathcal F}_2$ and
${\mathcal F}_1\cap{\mathcal F}_2$ are Cauchy filters on $X$.

\smallskip
Proof.} $(\Rightarrow)$ Suppose ${\mathcal F}_1$ and ${\mathcal
F}_2$ to be in the relation $R$ and fix $V\in{\mathcal U}$. Let
$W\in{\mathcal U}$ be such that $4W\subset V$. There exist
$F_1\in{\mathcal F}_1$ and $F_2\in{\mathcal F}_2$ such that
$F_1\times F_2\subset W$. Then $F_2\subset K(F_1,W)$ which implies
that $K(F_1,W)\in{\mathcal F}_2$. Since $F_1\subset K(F_1,W)$ we
have $K(F_1,W)\in{\mathcal F}_1$. Hence $K(F_1,W)\in{\mathcal
F}_1\cap{\mathcal F}_2$. On the other hand, for any $y_1,y_2\in
K(F_1,W)$ there are $x_1,x_2\in F_1$ such that
$(y_1,x_1),(y_2,x_2)\in W$. For an arbitrary chosen $z\in F_2$ we
have $(z,x_1),(z,x_2)\in W$. Hence $(y_1,y_2)\in 4W\subset V$. I
means that $K(F_1,W)\times K(F_1,W)\subset V$. Since $K(F_1,W)$ is
an element of ${\mathcal F}_1,{\mathcal F}_2$ and ${\mathcal
F}_1\cap{\mathcal F}_2$ we obtain all this filters to be Cauchy
filters.

$(\Leftarrow)$ Suppose ${\mathcal F}_1\cap{\mathcal F}_2$ to be
Cauchy filter on $X$. Fix $V\in{\mathcal U}$, choose
$F\in{\mathcal F}_1\cap {\mathcal F}_2$ such that $F\times
F\subset V$ and put $F_1:=F_2:=F$. Then $F_1\in {\mathcal F}_1,\
F_2\in {\mathcal F}_2$ and $F_1\times F_2\subset V$. Hence
${\mathcal F}_1R{\mathcal F}_2$. \hfill $\Box$\\

\smallskip
{\sc Proposition 3.4} {\it The relation $R$ described in
Definition 3.8 is an equivalence relation on the set $CF(X)$ of
all Cauchy filters on the uniform space $(X,{\mathcal U})$.

\smallskip
Proof.} It is obvious that for any Cauchy filters ${\mathcal F}_1$
and ${\mathcal F}_2$ on $X$ we have ${\mathcal F}_1R{\mathcal
F}_1$, and if ${\mathcal F}_1R{\mathcal F}_2$ then ${\mathcal
F}_2R{\mathcal F}_1$. Suppose now ${\mathcal F}_3$ to be such a
Cauchy filter on $X$ that ${\mathcal F}_1R{\mathcal F}_2$ and
${\mathcal F}_2R{\mathcal F}_3$. Fix $V\in{\mathcal U}$ and choose
$W\in{\mathcal U}$ such that $2W\subset V$. There exist
$F_1\in{\mathcal F}_1$, $F_2',F_2''\in{\mathcal F}_2$ and
$F_3\in{\mathcal F}_3$ such that $F_1\times F_2'\subset W$ and
$F_2''\times F_3\subset W$. Let $F_2:=F_2'\cap F_2''$. Then
$F_2\in{\mathcal F}_2$, $F_1\times F_2\subset W$ and $F_2\times
F_3\subset W$. Since $2W\subset V$ we obtain $F_1\times F_3\subset
V$. Hence ${\mathcal F}_1R{\mathcal F}_3$. \hfill $\Box$\\

\smallskip
For any Cauchy filter ${\mathcal F}$ on $X$ we denote by
$[{\mathcal F}]$ the equivalence class of ${\mathcal F}$ with
respect to the equivalence relation $R$ given in Definition 3.8.\\

\smallskip
{\sc Proposition 3.5} {\it If $\{{\mathcal F}_i\}_{i\in I}$ is a
family of Cauchy filters on an uniform space $X$ contained in an
equivalence class $[{\mathcal F}]$ then $\bigcap\limits_{i\in
I}{\mathcal F}_i\in[{\mathcal F}]$ i.e. $\bigcap\limits_{i\in
I}{\mathcal F}_i$ is a Cauchy filter and it is equivalent to
${\mathcal F}$.

\smallskip
Proof.} Let $V$ be an arbitrary element of the uniform structure
${\mathcal U}$ on $X$. Let $W\in{\mathcal U}$ be such that
$4W\subset V$. Choose $F\in{\mathcal F}$ such that $F\times
F\subset W$. Similarly as in the proof of Proposition 3.3 we
obtain that $K(F,W)\times K(F,W)\subset 3W\subset V$. For any
$i\in I$ there are $F_i\in{\mathcal F}_i$ and $G_i\in{\mathcal F}$
such that $F_i\times G_i\in W$. Hence $F_i\times(G_i\cap F)\in W$
and therefore $F_i\subset K(F,W)$. Consequently
$K(F,W)\in{\mathcal F}_i$ for any $i\in I$. Then
$K(F,W)\in\bigcap\limits_{i\in I}{\mathcal F}_i$ and moreover
$K(F,W)\times K(F,W)\subset V$. It means that
$\bigcap\limits_{i\in I}{\mathcal F}_i$ is a Cauchy filter on $X$.
Since $K(F,W)\in{\mathcal F}$ ($F\subset K(F,W)$) we obtain that
$\bigcap\limits_{i\in I}{\mathcal F}_i$ is equivalent to
${\mathcal F}$. \hfill $\Box$\\

\smallskip
{\sc Corollary 3.1} {\it If ${\mathcal F}$ is a Cauchy filter on
$X$ then $\bigcap\limits_{{\mathcal G}\in [{\mathcal F}]}{\mathcal
G}$ is a Cauchy filter on $X$ equivalent to ${\mathcal F}$. Since
for any ${\mathcal F}_1 \in [{\mathcal F}]$ we have
$\bigcap\limits_{{\mathcal G} \in [{\mathcal F}]} {\mathcal G}
\subset {\mathcal F}_1$ we obtain $\bigcap\limits_{{\mathcal G}
\in [{\mathcal F}]} {\mathcal G}$ is the minimal element of $[{\mathcal F}]$
with respect to the ordering relation $\subset$ on the family of
all filters on a set $X$.} \hfill $\Box$\\

\smallskip
{\sc Definition 3.9} For any Cauchy filter ${\mathcal F}$ on $X$
the Cauchy filter $\bigcap\limits_{{\mathcal G}\in [{\mathcal
F}]}{\mathcal G}$ is called {\sl the minimal Cauchy filter on $X$
defined by (smaller then)} ${\mathcal F}$.

\bigskip
{\sc Definition 3.10}  A uniform space $(X,{\mathcal U})$
is said to be {\sl complete} if each Cauchy filter on $X$ is
convergent in $\tau_{\mathcal U}$.

\bigskip
{\sc Theorem 3.2} {\it If $(X,{\mathcal U})$ is a complete uniform
space and $M$ is a closed subset of the topological space
$(X,\tau_{\mathcal U})$ then a uniform space $(M,{\mathcal U}_M)$
is complete. Conversely, if $(M,{\mathcal U}_M)$ is a complete
uniform subspace of some (not necessarily complete) uniform space
$(X,{\mathcal U})$, then M is closed in $X$ with respect to}
$\tau_{\mathcal U}$.

\smallskip
For the proof see \cite{Bour}, \cite{eng} or \cite{pas}.\\

\smallskip
It is well known that the uniform space of reals $({\mathbb
R},{\cal U})$ with the standard uniform structure ${\cal U}={\cal
U}_{\{id_{\mathbb R}\}}$ defined by the one element family of
functions $\{id_{\mathbb R}\}$ (or by the standard metric) is
complete. We have also more general

\smallskip
{\sc Proposition 3.6} {\it For any set $I$ the uniform space
$({\mathbb R}^I,{\cal U}_{\cal G})$, where ${\cal
G}=\{pr_i\}_{i\in I}$ is the set of all natural projections
$pr_i:{\mathbb R}^I\to{\mathbb R}$,
$$pr_i(f)=f(i),\qquad f\in{\mathbb
R}^I,$$ for any $i\in I$, is complete.

\smallskip
Proof.} If ${\cal F}$ is a Cauchy filter on ${\mathbb R}^I$ then
for any $i\in I$ the set $pr_i({\cal F})=\{pr_i(F):F\in{\cal F}\}$
is a filtering base of some Cauchy filter on ${\mathbb R}$. Then
the Cauchy filter corresponding to $pr_i({\cal F})$ converges to
some $y_i\in{\mathbb R}$. Putting $f(i):=y_i,\ \ i\in I$ we obtain
function $f\in{\mathbb R}^I$ such that ${\cal F}\rightarrow f$.
\hfill$\Box$

\bigskip
Any uniform space can be treated as a uniform subspace of some
complete uniform space. We have the following

\smallskip
{\sc Theorem 3.3} {\it  For each uniform space} $(X,{\mathcal U})$:\\
(i) {\it there exists a complete uniform space
$(\widetilde{X},\widetilde{\mathcal U})$ and a set
$A\subset\widetilde{X}$ dense in $\widetilde{X}$ (with respect to
the topology $\tau_{\widetilde{\mathcal U}}$) such that
$(X,\mathcal{U})$ is uniformly homeomorphic to}
$(A,\widetilde{\mathcal U}_A)$;\\
(ii) {\it if the complete uniform spaces
$(\widetilde{X}_1,\widetilde{\mathcal U}_1)$ and
$(\widetilde{X}_2,\widetilde{\mathcal U}_2)$ satisfies condition
of the point} (i) {\it then they are uniformly homeomorphic}.

\smallskip
For the details of the proof see \cite{Bour} or \cite{pas}. Here
we only want to describe the construction of
$(\widetilde{X},\widetilde{\mathcal U})$.

\bigskip
Let $\widetilde{X}$ be the set of all minimal Cauchy filters in
$X$. For every $V\in{\mathcal U}$ we denote by $\widetilde{V}$ the
set of all pairs $({\mathcal F}_1,{\mathcal F}_2)$ of minimal
Cauchy's filters, which have a common element being a small set of
rank V. We define a family $\widetilde{\mathcal U}$ of subsets of
set $\widetilde{X}\times \widetilde{X}$ as the smallest uniform
structure on $X$ containing all sets from the family
$\{\widetilde{V}:V\in{\mathcal U}\}$.

\bigskip
{\sc Example 3.3} Let us consider two  uniform structures
${\mathcal U}_{\{f\}}$ and ${\mathcal U}_{\{g\}}$ on the
differential space $({\mathbb R},{\mathcal C}^{\infty})$, where
$$f(x)=x,\ \ g(x)=arctg x,\qquad x\in{\mathbb R}.$$
Then $({\mathbb R},{\mathcal U}_{\{f\}})$ is the complete space
i.e.  $\widetilde{{\mathcal U}_{\{f\}}} = {\mathcal U}_{\{f\}}$
while $({\mathbb R}, {\mathcal U}_{\{g\}})$ is not complete and we
have $\widetilde{X} \simeq[- \frac{\pi}{2} ; \frac{\pi}{2} ]$.
Consequently ${\mathcal U}_{\{f\}}\neq{\mathcal U}_{\{g\}}$.

\bigskip
Let $N$ be a set, $M \subseteq N$, $M \neq \O$, ${\mathcal C}$ be
a differential structure on $M$.

\smallskip
{\sc Definition 3.11.} The differential structure ${\mathcal D}$
on $N$ is an \emph{extension} of the differential structure
${\mathcal C}$ from the set $M$ to the set $N$ if ${\mathcal C} =
{\mathcal D}_{M}$ (if we get the structure ${\mathcal C}$ by
localization of the structure ${\mathcal D}$ to $M$).

\bigskip
For the sets $N,M$ and the differential structure ${\mathcal C}$
on $M$ we can construct many different extensions of the structure
$M$ to $N$.

\smallskip
{\sc Example 3.4.} If for each function $f\in {\mathcal C}$ we
assign $f_{0} \in {\mathbb R}^{N}$ such that $f_{0|M} = f$ and
$f_{0|N \backslash M} \equiv 0$, then the differential structure
generated on $N$ by the family of functions $\{f_{0}\}_{f \in
{\mathcal C}}$ is an extension of ${\mathcal C}$ from $M$ to $N$.
Similarly, if for each function $f \in {\mathcal C}$ we assign the
family ${\mathcal F}_{f} := \{g \in {\mathbb R}^{N} : g_{|M} =
f\}$, then the differential structure on $N$ generated the family
of functions ${\mathcal F} := \bigcup_{f \in {\mathcal C}}
{\mathcal F}_{f}$ is an extension of ${\mathcal C}$ from $M$ to
$N$. If the set $N \backslash M$ contain at least two elements,
then the differential structures generated by the families
$\{f_{0}\}_{f \in {\mathcal C}}$ and ${\mathcal F}$ are different.

\bigskip
{\sc Definition 3.12.} If $\tau$ is a topology on the set $N$,
then the extension ${\mathcal D}$ of the differential structure
${\mathcal C}$ from $M$ to $N$ is \emph{continuous with respect to
$\tau$} if each function $f \in {\mathcal D}$ is continuous with
respect to $\tau$ $(\tau_{\mathcal D} \subset\tau)$.

\bigskip
If on the set $N$ there exists a continuous with respect to some
topology $\tau$ extension of the differential structure ${\mathcal
C}$ from the set $M \subset N$, then the structure ${\mathcal C}$
is said to be {\sl extendable from the set $M$ to the topological
space} $(N, \tau)$.

\bigskip
{\sc Example 3.5.} The differential structure $C^{\infty}({\mathbb
R})_{\mathbb Q}$ is extendable from the set of rationales to the
set of reals. The continuous extensions are e.g.
$C^{\infty}({\mathbb R})$ and the structure ${\mathcal D}$
generated on ${\mathbb R}$ by the family of the functions
$C^{\infty}({\mathbb R}) \cup\{f\}$, where $f : {\mathbb R}
\rightarrow {\mathbb R}$, $f(x) := |x - \sqrt{2}|$, $x \in
{\mathbb R}$.

\bigskip
{\sc Proposition 3.7} {\it Let $M \neq \emptyset$, $(M, \mathcal
C)$ be a differential space and ${\cal G}$ be a family of
generators of ${\mathcal C}$, i.e. ${\mathcal C} = gen({\cal G})$.
Let $(\widetilde{M},\widetilde{\mathcal U}_{\cal G})$ be the
completion of the uniform space $(M,{\cal U}_{\cal G})$. Then any
function $g\in{\cal G}$ poses  the continuous extension
$\tilde{g}:\widetilde{M}\to{\mathbb R}$. If $\widetilde{\cal G}$
is the family of all continuous extensions of elements of ${\cal
G}$ to $\widetilde{M}$ then the differential structure ${\cal
D}=gen(\widetilde{\cal G})$ is an continuous extension of the
differential structure ${\mathcal C}$ from the set $M$ to the set
$\widetilde{M}$. Moreover $\tau_{\cal D}=\tau_{\tilde{\cal
U}_{\cal G}}$.

\smallskip
Proof.} Let $\phi_{{\mathcal G}}$ be the generator embedding of
the differential space $(M,{\mathcal C})$ into the Cartesian space
$({\mathbb R}^{\cal G},C^\infty({\mathbb R}^{\cal G}))$ defined by
${\mathcal G}$.  Then the closure $\overline{\phi_{{\mathcal
G}}(M)}$ is a complete subspace of the complete uniform space
${\mathbb R}^{\cal G}$ (see Theorem 3.2 and Proposition 3.6). We
know that ${\mathbb R}^{\cal G}:(M,{\mathcal
C})\to(\phi_{{\mathcal G}}(M),C^\infty({\mathbb R}^{\cal
G})_{\phi_{{\mathcal G}}(M)})$ is a diffeomorphism and
$pr_g\circ\phi_{{\mathcal G}}=g$ for any $g\in{\cal G}$. Moreover
$\phi_{{\mathcal G}}(M)$ is dense in $\overline{\phi_{{\mathcal
G}}(M)}$. Then identifying any $g\in{\cal G}$ with
${\pi_g}_{|\phi_{{\mathcal G}}(M)}$, ${\cal C}$ with
$C^\infty({\mathbb R}^{\cal G})_{\phi_{{\mathcal G}}(M)}$ and
putting $\widetilde{M}:=\overline{\phi_{{\mathcal G}}(M)}$ we
obtain that $\tilde{g}$ should be identify with
${pr_g}_{|\overline{\phi_{{\mathcal G}}(M)}}$ and ${\cal
D}=C^\infty({\mathbb R}^{\cal G})_{\overline{\phi_{{\mathcal
G}}(M)}}$. We have also $\tau_{\cal
D}=\tau_{\widetilde{M}}=\tau_{\widetilde{\cal U}_{\cal G}}$, where
$\tau_{\widetilde{M}}$ is the topology of $\widetilde{M}$ treated
as a topological subspace of ${\mathbb R}^{\cal G}$. \hfill $\Box$

\bigskip
{\sc Definition 3.12.} The differential space
$(\widetilde{M},{\cal D})$ constructed in Proposition 3.7 will be
called {\sl the differential completion} of the differential space
$(M,{\cal C})$ defined by the family of generators ${\cal G}$. The
set $\widetilde{M}$ will be denoted by $compl_{\cal G}M$ and the
differential structure ${\cal D}$ will be denoted by $compl_{\cal
G}{\cal C}$.

\section {The maximal differential completion}

Let us consider two families ${\cal G}$ and ${\cal H}$ of
generators of a differential structure ${\cal C}$ on a set
$M\neq\emptyset$. If ${\cal G}\subset{\cal H}$ then for uniform
structures ${\cal U}_{\cal G}$ and ${\cal U}_{\cal H}$ we have:
${\cal U}_{\cal G}\subset{\cal U}_{\cal H}$. Consequently any
Cauchy filter with respect to ${\cal U}_{\cal H}$ is a Cauchy
filter with respect to ${\cal U}_{\cal G}$. In particular any
minimal Cauchy filter with respect to ${\mathcal U}_{\mathcal H}$
is a Cauchy (but not necessarily minimal Cauchy) filter with respect
to ${\mathcal U}_{\mathcal G}$. This defines the natural map
$\iota_{{\mathcal G}{\mathcal H}} : compl_{\mathcal H} M \rightarrow compl_{\mathcal G} M$
as follows: for any ${\mathcal F} \in compl_{\mathcal H} M$ the value
$\iota_{{\mathcal G}{\mathcal H}} ({\mathcal F}) \in compl_{\mathcal G} M$
is the minimal Cauchy filter equivalent to ${\mathcal F}$ with respect
to the uniform structure $\mathcal{U}_{\mathcal{G}}$.

\bigskip
{\sc Proposition 4.1} {\it For any two families $\mathcal{G}$ and $\mathcal{H}$
of generators of a differential structure $\mathcal{C}$ on a set
$M \neq \emptyset$ such that $\mathcal{G} \subset \mathcal{H}$ the map
$\iota_{{\mathcal G}{\mathcal H}} : compl_{\mathcal H} M \rightarrow compl_{\mathcal G} M$
defined above is smooth with respect to differential structures
$compl_{\mathcal{H}}{\mathcal{C}}$ and $compl_{\mathcal{G}}{\mathcal{C}}$.

\smallskip
Proof.} For smoothness of $\iota_{{\mathcal G}{\mathcal H}}$ it is
enough to prove that for any $g \in \mathcal{G}$ the function
$\tilde{g}_{\mathcal{G}} \circ \iota_{{\mathcal G}{\mathcal H}}
\in compl_{\mathcal{H}}{\mathcal{C}}$, where
$\tilde{g}_{\mathcal{G}}$ denotes the continuous extension of $g$
onto $compl_{\mathcal{G}}M$. Since any $g \in{\mathcal G}$ is an
element of $\mathcal{H}$ we have for each Cauchy filter ${\mathcal
F} \in compl_{\mathcal H}M$
$$\tilde{g}_{\mathcal G} \circ \iota_{{\mathcal G}{\mathcal H}}({\mathcal F}) =
\lim{g(\iota_{{\mathcal G}{\mathcal H}}({\mathcal F}))} = \lim{g({\mathcal F})}
= \tilde{g}_{\mathcal H}({\mathcal F}),$$
where $\tilde{g}_{\mathcal{H}}$ denotes the continuous extension of $g$ onto
$compl_{\mathcal{H}}M$. Hence $\tilde{g}_{\mathcal{G}} \circ \iota_{{\mathcal G}
{\mathcal H}} = \tilde{g}_{\mathcal{H}} \in compl_{\mathcal{H}}{\mathcal{C}}$.

\bigskip
{\sc Remark 4.1} {\it In general the image $\iota_{{\mathcal G}{\mathcal H}}
(compl_{\mathcal{H}}M) \neq compl_{\mathcal{G}}M$ and it is not complete in
$compl_{\mathcal{G}}M$. For example let $M = (0;\frac{\pi}{2})$, ${\mathcal C}
= C^\infty (M)$ and ${\mathcal G} = \{x \cos(\tan x), x \sin (\tan x)\}$,
${\mathcal H} = \{x \cos (\tan x), x \sin (\tan x), \tan x\}$.}

\bigskip
Then we obtain the following

\bigskip
{\sc Theorem 4.1} {\it For any differential space $(M,{\cal C})$ the differential
completion $(compl_{\cal C}M,compl_{\cal C}{\cal C})$ has the following properties:\\
(i) for any differential completion $(N,{\cal D})$ of  $(M,{\cal
C})$ (where $M \subset N$) there exist the map
$$\iota_{\mathcal{D}} : (compl_{\mathcal C}M, compl_{\mathcal C}{\mathcal C}) \rightarrow (N,\mathcal{D})$$
such that $\iota_{\mathcal{D}|M} = id_M$;\\
(ii) for any function $g\in{\cal C}$ there exists uniquely defined extension $\tilde{g}\in compl_{\cal C}{\cal C}$.}

\bigskip
In the set of all differential completions of the space $(M,\mathcal{C})$ we can define
an ordering relation $\preceq$ such that:
$$compl_{\mathcal G}M \preceq compl_{\mathcal H}M \Leftrightarrow {\mathcal G} \subset {\mathcal H},$$
where ${\mathcal G}$ and ${\mathcal H}$ are families of generators of the structure ${\mathcal C}$.

\smallskip
The above theorem says that $compl_{\cal C}M$ is the maximal with the respect to the order
$\preceq$ completion of $M$ which
can be constructed using a set of generators of the differential structure ${\cal C}$
while $compl_{\cal C}{\cal C}$ is the maximal continuous extension of ${\cal C}$ from
$M$ to $compl_{\cal C}M$.

\bigskip
{\sc Definition 4.1.} We will call the differential space $(compl_{\cal C}M,compl_{\cal
C}{\cal C})$ {\sl the maximal differential completion} of the differential space
$(M,{\cal C})$.

\bigskip
Let us consider the situation when for some family of generators ${\cal G}$ the
uniform space $(M,{\cal U}_{\cal G})$ is complete.

\smallskip
{\sc Theorem 4.2} {\it Let $(M, \mathcal C)$ be a differential space and ${\cal G}$ be
a family of generators of ${\mathcal C}$. If the uniform space $(M,{\cal U}_{\cal G})$
is complete then for any family ${\cal H}$ of generators of ${\mathcal C}$ such that
${\cal G}\subset {\cal H}$ we have
\begin{equation} \label{MAX1}
compl_{\cal H}M=M
\end{equation}
and
\begin{equation} \label{MAX2}
compl_{\cal H}{\cal C}={\cal C}.
\end{equation}
In particular $compl_{\cal C}M=M$ and $compl_{\cal C}{\cal
C}={\cal C}$.

\smallskip
Proof.} Any element of $compl_{\cal H}M$ is represented by some filter ${\cal F}$ in $M$
which is a Cauchy filter with respect to ${\cal U}_{\cal H}$. Then ${\cal F}$ is a Cauchy
filter with respect to ${\cal U}_{\cal G}$ and therefore ${\cal F}$ can be identify
with its limit in M. Then $compl_{\cal H}M\subset M$ On the other hand for any element
$p\in M$ the filter ${\cal F}_p$ of all neighbourhoods of $p$ is a Cauchy filter with
respect to ${\cal U}_{\cal H}$. Hence we can write $M\subset compl_{\cal H}M$.

The equality (\ref{MAX2}) is an immediate consequence of the definition of
$compl_{\cal H}{\cal C}$ and the equality (\ref{MAX1}).  \hfill $\Box$

\bigskip
Let us consider the case when ${\cal G}={\cal C}$.

\smallskip
{\sc Theorem 4.3} {\it Let $(M,\mathcal C)$ be a differential space If there exists
a finite or countable family of generators of ${\cal C}$ then the uniform space
$(M,{\cal U}_{\cal C})$ is complete.

\smallskip
Proof.} Suppose that $(M,{\cal U}_{\cal C})$ is not complete. Then
there exists $x\in compl_{\cal C}M\setminus M$. Let $\chi:{\cal
C}\to{\mathbb R}$ be a functional given by the formula
\begin{equation} \label{MAX3}
\chi(g):=\tilde{g}(x), \qquad g\in{\cal C},
\end{equation}
where $\tilde{g}$ is the continuous extension of $g$ from $M$ to
$compl_{\cal C}M$. This functional is an element of the spectrum
of the algebra ${\cal C}$ but it is not an evaluation functional
on $M$ (the algebra $compl_{\cal C}{\cal C}$ separates points of
the space $compl_{\cal C}M$). Then ${\cal C}$ does not posses the
spectral property. It is a contradiction with Theorem 1 and
Corollary 6 from the work \cite{cukPWS} (see also Theorem 2.3
(Twierdzenie 2.3) and Corollary 2.6 (Wniosek 2.6) from
\cite{cuk}). \hfill $\Box$

\bigskip
{\sc Corollary 4.1} {\it Let $X$ be a topological Hausdorff space.
If the topology of $X$ is given by a countable family ${\cal G}$
of real-valued functions as the weakest topology on $X$ with
respect to which all elements of ${\cal G}$ are continuous then
there exists an uniform structure ${\cal U}$ on $X$ such that the
uniform space $(X,{\cal U})$ is complete and the topology
$\tau_{\cal U}$ coincides with the initial topology on $X$.}
\hfill $\Box$

\section{Compactification of a differential space}

Let $(M, {\mathcal C})$ be a differential space such that
${\mathcal C} = gen({\mathcal G})$. Let $f \in {\mathcal C}$, $m
\in M$. Then there exist: neighbourhood $U$ of $m$, number $n
\in{\mathbb N}$, functions $\alpha_1, \ldots, \alpha_n \in
{\mathcal G}$ and $\omega \in C^{\infty}({\mathbb R}^n)$ such that
$f|_{U} = \omega \circ (\alpha_1, \ldots, \alpha_n)|_{U}$. We
denote $y_0:= (\alpha_1(m), \ldots, \alpha_n(m))$. Let us take cubes:\\
$P:=(\alpha_1(m)-1, \alpha_1(m)+1) \times (\alpha_2(m)-1,
\alpha_2(m)+1) \times \ldots \times (\alpha_n(m)-1, \alpha_n(m)+1),$\\
$P':=(\alpha_1(m)-2, \alpha_1(m)+2) \times (\alpha_2(m)-2,
\alpha_2(m)+2) \times \ldots \times (\alpha_n(m)-2,
\alpha_n(m)+2)$.

Let $\eta \in C^\infty ({\mathbb R}^n)$, such that: $\eta|_P
\equiv 1$, $\eta|_{{\mathbb R}^{n} \setminus P'} \equiv 0$ and
$|\eta| \leq 1$. We mark: $\alpha:=(\alpha_1, \ldots , \alpha_n)$,
$\beta_i:=\alpha_i \cdot \eta (\alpha_1, \ldots , \alpha_n)$,
$\beta:=(\beta_1, \ldots, \beta_n)$. Let $V:=\alpha^{-1}(P)$, $V'
:=\alpha^{-1}(P')$. Then $V \subset V' \subset M$ and $m \in U
\cap V$. For any $x \in U \cap V$ we have: $f(x) = \omega
(\alpha_{1}(x), \ldots, \alpha_{n}(x)) = \omega (\beta_{1}(x),
\ldots, \beta_{n}(x))$. We observe that $\forall i \in \{1,
\ldots, n\}$ $|\beta_{i}(x)| \leq max\{|\alpha_{1}(m) + 2|,
|\alpha_{1}(m) - 2|\} =: \mu_i$. Then $f(x) = \omega (\mu_1
\frac{\beta_1 (x)}{\mu_1}, \mu_2 \frac{\beta_2 (x)}{\mu_2},
\ldots, \mu_n \frac{\beta_n (x)}{\mu_n}) = \omega (\mu_1 \gamma_1
(x), \ldots, \mu_n \gamma_n (x)) = \omega_1 (\gamma_1, \ldots,
\gamma_n)(x)$, where: $\forall 1 \leq i \leq n$, $\forall x \in M$
$\gamma_i (x) = \frac{\beta_i (x)}{\mu_i}$ and $|\gamma_i (x)|<
1|$. Hence we get the following theorem:

\bigskip
{\sc Theorem 5.1} {\it For any differential space $(M, {\mathcal
C})$ there exist the family of bounded generators, in particular
${\mathcal C} = gen({\mathcal G}_{1})$, where ${\mathcal G}_{1} =
\{\gamma_i\}_{|i \in I}$ such that $|\gamma_i| \leq 1$ $\forall i
\in I$.}

\bigskip
For any Hausdorff differential space we consider the generators
embedding of that space using the family of generators described
in Theorem 5.1 (it takes values in the cube $[-1, 1]^I$). So we
have the embedding of the differential space into the compact
space and we close the image. Hence we get the compact set
$\overline{M}_{\mathcal G}$.

If we mark $J:=[-1, 1]$, then $\overline{M}_{\mathcal G} \subset
J^I$. On $J^I$ there exists the natural differential structure
$C^\infty (J^I) = C^\infty ({\mathbb R}^I)|_{J^I}$ generated by
the family of projections $\{pr_{i|_J}\}_{i \in I}$, where $pr_i :
{\mathbb R}^I \rightarrow{\mathbb R}$ is the projection onto
$i$-th coordinate. By localization of that structure to the set
$\overline{M}_{\mathcal G}$ we get the differential space
$(\overline{M}_{\mathcal G}, C^\infty
(J^I)_{\overline{M}_{\mathcal G}})$ which is a differential
subspace of the space $(J^I, C^\infty (J^I))$. We call that
differential space \emph{the (differential) compactification of
the differential space $(M, {\mathcal C})$ by the family of
generators} ${\mathcal G}$ and we denote it by $(compt_{\mathcal
G}M, compt_{\mathcal G}{\mathcal C})$. We see that
$compt_{\mathcal G}M = compl_{\mathcal G}M$ and $compt_{\mathcal
G}{\mathcal C} = compl_{\mathcal G}{\mathcal C}$. So for the
compactification of the differential space we have analogous
theorems like for the completion.

\bigskip
{\sc Remark 5.1} {\it Let ${\mathcal G}$ and ${\mathcal H}$ be the
families of bounded generators of the differential structure
${\mathcal C}$ on the set $M \neq \O$. If ${\mathcal G} \subset
{\mathcal H}$, that there exists smooth function $\iota_{{\mathcal
G}{\mathcal H}} : (compt_{\mathcal H}M, compt_{\mathcal
H}{\mathcal C}) \rightarrow (compt_{\mathcal G}M, compt_{\mathcal
G}{\mathcal C})$ such that $\iota_{{\mathcal G}{\mathcal H}|M} =
id_M$.}

\bigskip
Let consider the differential space $(M, {\mathcal C})$ and the
family ${\mathcal C}_0$ of all smooth functions on that space
which takes values in $[ -1, 1 ]$. Using the procedure of the
compactification, described earlier, we get the differential space
that we mark by $(compt_{\mathcal C}M, compt_{\mathcal C}{\mathcal
C})$ it means $compt_{\mathcal C}M := compt_{{\mathcal C}_0}M$,
$compt_{\mathcal C}{\mathcal C} := compt_{{\mathcal C}_0}{\mathcal
C}$.

\bigskip
{\sc Definition 5.1} The differential space $(compt_{\mathcal C}M,
compt_{\mathcal C}{\mathcal C})$ is called the \emph{maximal
differential compactification} of the space $(M, {\mathcal C})$

\bigskip
Let us assume that the topological space $(M, \tau_{\mathcal C})$
is compact. Than all the functions from ${\mathcal C}$ are limited
and by the normalization of each function $\alpha \in {\mathcal C}
\setminus \{0\}$ according the formula:
$$N {\alpha(p)} = \frac{1}{sup_{q \in M}|\alpha(q)|} \alpha(p), p \in M$$
we get the family $N{\mathcal C} = \{N \alpha : \alpha \in
{\mathcal C}\}$ of the generators of the structure ${\mathcal C}$.
Then the generators embedding given by $N{\mathcal C}$ converts
diffeomorphically $(M, {\mathcal C})$ onto the compact subspace of
$(J^{N{\mathcal C}}, C^\infty (J^{N{\mathcal C}}))$. Similarly, if
${\mathcal G}$ is any family of generators of the structure
${\mathcal C}$, then $N{\mathcal G} = \{N \alpha : \alpha \in
{\mathcal G}\}$ is the family of generators of ${\mathcal C}$ too,
and proper generators embedding is the diffeomorphism of $(M,
{\mathcal C})$ onto a compact subspace of $(J^{N{\mathcal G}},
C^\infty (J^{N{\mathcal G}}))$. We have
$$compt_{\mathcal C}M = M, \quad compt_{\mathcal G}{\mathcal C} =
{\mathcal C},$$ where the equality is the identification of the
diffeomorphic spaces and structures.

 \vspace{1cm}

\footnotesize
\begin{center}

\end{center}


\begin{thebibliography}{99}

\bibitem{Bour} N. Bourbaki, {\sl General Topology}, Springer-Verlag,
Berlin Heidelberg New York London Paris Tokyo 1989.

\bibitem{cuk} M. J. Cukrowski, {\sl Geometrical properties of function
algebras in the category of differential spaces (doctor thesis -
in polish)}, Warsaw University of Technology, Warsaw 2007.

\bibitem{cukPWS} M. J. Cukrowski, Z. Pasternak-Winiarski, W. Sasin
{\sl On real-valued homomorphisms in countably generated
differential structures}, arXiv:1103.3597, 2011.

\bibitem{DPW2} D. Dziewa-Dawidczyk, Z. Pasternak - Winiarski, {\sl Differential
structures on natural bundles connected with a differential
space}, in "Singularities and Symplectic Geometry VII" Singularity
Theory Seminar (2009), S. Janeczko (ed), Faculty of Mathematics
and Information Science, Warsaw University of Technology

\bibitem{DPW1} D. Dziewa-Dawidczyk, Z. Pasternak - Winiarski, {\sl Uniform structures
on differential spaces}, arXiv:1103.2799, 2011.

\bibitem{eng} R. Engelking, {\sl General topology}, PWN, Warsaw
1975  (in polish).

\bibitem{pas} Z. Pasternak - Winiarski, {\sl Group differential structures
and their basic properties (doctor thesis - in polish)}, Warsaw
University of Technology, Warsaw 1981.

\bibitem{sik} R. Sikorski, {\sl Introduction to the differential geometry}, PWN, Warsaw
1972 (in polish).

\bibitem{wal} W. Waliszewski, {\sl Regular and coregular mappings of differential
space}, Annales Polonici Mathematici XXX, 1975.

\end{thebibliography}
\end{document}